# Analysis of the complex gas pipeline exploitation process in various operating modes.

**Ilgar Aliyev**

**Head of Operation and reconstruction of buildings and facilities Department, Azerbaijan University Architecture and Construction, Baku, Azerbaijan**

**Abstract:** The study aims to decrease gas loss and enhance system reliability during gas pipeline accidents. A computational scheme has been developed that can enable the elimination of gas leakage through the modeling and management of parallel gas pipeline systems. The dynamic state of processes for the supply of modern automatic equipment to gas pipelines and the use of an efficient automated control system have been extensively studied. The analytical determination of the optimal transition time has been widely applied to ensure the most favorable operating conditions for the system. Methods for calculating complex transient processes in main gas pipelines, from a non-stationary regime to a stationary regime, have been developed, particularly at the moment of gas flow ingress. A comparison of mathematical expressions for calculating transient processes in complex main gas pipelines has been conducted through theoretical sources.

**Introduction:** The accidents occurring in the exploited main gas pipelines cause significant losses to the country's economy and lead to environmental pollution. Therefore, the development and preparation of effective methods for the analysis of non-stationary regimes of natural gas transportation in main gas pipelines undoubtedly arouses scientific and practical interest. At the same time, solving dynamic problems in complex main gas pipelines will allow discovering previously unknown relationships and legalities, conducting a more thorough analysis of the problem, and providing substantial results and recommendations. The relevance of this problem has increased in recent years, especially in connection with the ever-increasing requirements for the accuracy of calculations for automated control systems of technological processes. Calculation of transient processes in the main gas pipeline section is necessary for determining the project and exploitation time pressure limit and for developing and preparing automatic protection means in case of exceeding this limit [1,2]

Indicators of pressure increase or decrease in the compressor station become known at the moment of activating the automatic protection system and regulating it until the release of protective devices necessary for protection. Furthermore, studying transient processes in complex gas pipelines during complex situations, such as pressure increase or decrease, is crucial for identifying and resolving accidents in a timely manner in the main gas pipelines. This is especially important for the operational dispatch control of automated technological management systems in the automated technological management systems of main gas pipelines [3].

However, one of the challenges in calculating complex transient processes in main gas pipelines during the non-stationary regime, when gas flow transitions from a non-stationary regime to a stationary regime in the pipeline (according to the proposed scheme), is the diverse parameters of the interconnected sections of the pipelines at the moment of gas flow transition [4].The nature of the challenges lies in the excessive complexity of the mathematical representation of dynamic processes in these systems. The difficulty arises because the distance between connectors or operating valves is identified after the location of the leakage is known, and these processes are described by heat conduction equations. In other words, all processes are described by non-linear differential equations. Computational methods are available for determining the lengths of the steps between connectors, and they are installed depending on the length of the main gas pipeline and the designation of gas consumers [5]. However, the quantity of gas flow transferred from the damaged section to the undamaged section of the pipeline and the law of change are not known in advance and should be determined according to the location of the gas pipeline leakage. Therefore, the methods for calculating transient processes in complex main gas pipelines should be sufficiently developed to allow for both their design and utilization processes, as well as the development and preparation of automatic control means.

**Materials, Methods and Results**

In the practice of engineering calculations, the analysis of transient processes in complex main gas pipelines is frequently encountered when the condition of installed valves at various points changes (when opening or closing) [7]. The reason for re-examining this issue is to address a series of problems for the automated control systems of technological processes in pipeline chambers, which result from determining the course of the technological process and the sequence of operation of technological devices. Considering this point, when connecting a parallel gas pipeline to exploitation, the calculation method should be applied to determine the sequence of operation (opening or closing) of valves installed on connectors and pipelines, as well as to determine the pressure value at characteristic points of the parallel main gas pipeline. For this purpose, a new calculation scheme has been proposed for the design and reconstruction of the parallel main gas pipeline that perfectly accommodates automated control systems for technological processes (Fig.1).

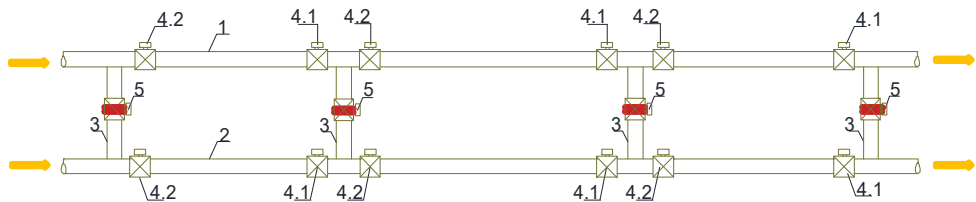

**Figure 1.** *Proposed schematic diagram of the parallel main gas pipeline.*

*1,2 - pipeline chambers on level 1 and level 2; 3 - connectors; 4.1 and 4.2 - automatic valves on pipelines; 5 - automatic valves on connectors.*

It is clear that fittings (valves, drawers, automatic taps, control valves, etc.) are installed on gas pipelines in order to control the flow of gas transported through pipelines, to separate one part of the pipeline from another, and to ensure the operation modes of technological equipment.

The technological aspects of the proposed reporting scheme that distinguish from the scheme of existing complex gas pipelines are as follows:

1. The distance (step) between pipes connecting the parallel pipeline networks of the complex gas pipeline is determined in advance, and the automatic valves installed on these connectors are in a closed position in the normal operating mode (stationary mode) of the gas pipeline. Only when there is a disruption in the integrity of the pipeline networks or during repairs are the automatic valves installed on the connectors activated, i.e., opened, and after the gas pipeline returns to normal operating mode, they are closed again.
2. In the pipeline networks, automatic valves are installed on the right side after the first connector and on the left side before the last connector. It is proposed to install valves on both the left and right sides of the other connectors (Scheme 1). Thus, in emergency modes, the activation of valves installed on connectors and pipeline networks creates a loop system, resulting in uninterrupted gas supply to consumers from the gas pipeline (Fig. 2).

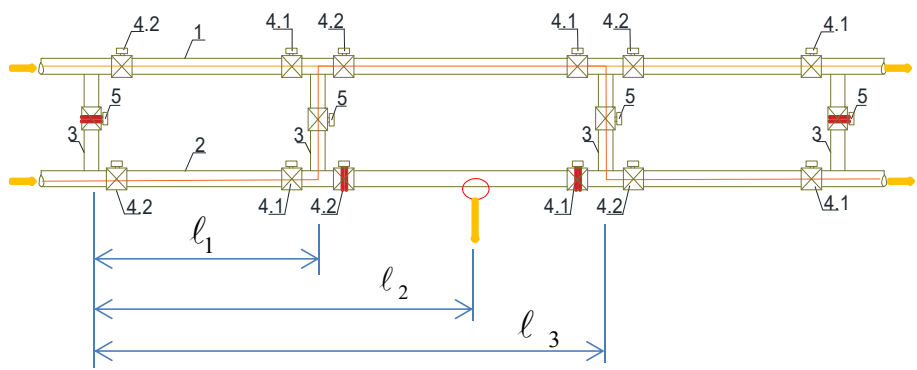

**Figure 2.** *Principal scheme of the proposed technological operating mode for the complex gas pipeline in emergency mode.*

From Figure 2, it is evident that in the event of a breach at point x=$\ell_2$ in the gas pipeline, the sequential operation of the technological operating mode for automated control systems of technological processes is implemented as follows. The scheme shows that two automatic valves are placed on each side of a connector, both on the left and right sides. However, in the event of an emergency, only one of these valves should operate. Therefore, the operational principle of these valves should be controlled not based on pressure drop but through the application of automated control systems in the emergency dispatching centers.

To achieve this, during an emergency, the location of the automatic valves that will operate in the control center must be known first. Initially, the location of the pipeline integrity breach (leakage point x=$\ell_2$) is identified in the emergency dispatching center. As mentioned, the length of the step between the valves becomes known from the moment the gas pipeline is put into operation. Assuming we know this distance, the following calculation method should be provided to determine the location of the valves ($\ell_3$ and $\ell_1$) installed on the left -4.2 and right -4.1 parts of point x=$\ell_2$ in the control center. To first determine the location of the valve situated on the left side of the emergency point, the following statement is designated:.

$$\ell_3 = \frac{\ell_2}{\ell} + a \quad \ell, m$$

Here, $a \in \,]0:1[$, İn other words $\frac{\ell_2}{\ell} + a$ =n: n=1,2,3,4......

Then, $\ell_1 = \ell_3 - \ell$

Thus, as the leakage point is identified in the emergency dispatch center, the locations of the automatic valves installed on the right and left sides of the damaged pipeline are determined and closed. As a result, the damaged (integrity compromised) section of the pipeline is separated from the main section of the gas pipeline. As a result, the valves on the connectors with known positions next to the 4.2 and 4.1 automatic valves are sequentially opened, allowing the gas flow from the damaged section of the gas pipeline to be directed to the undamaged section. This ensures uninterrupted gas supply to consumers relying on the gas pipeline.

After the closure of the 4.2 and 4.1 automatic valves, the transition process in the gas pipeline is divided into three sections. In the first section ($0 \le x \le \ell_1$), a filling process occurs in the pipeline, and the pressure value along the pipeline increases over time. The pressure increase affects the compression ratio of the gas compressors at the compressor station and can lead to a disruption of their operating mode (Fig. 3).

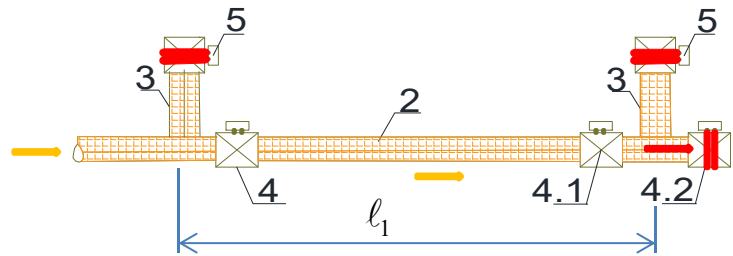

**Figure 3.** Principal scheme of the technological operating mode in the first section of the complex gas pipeline during the accident regime.

In order to prevent the progression of this process, it is essential to channel the compressed gas flow to the undamaged pipeline. For this purpose, the displacements installed on the connectors should be activated (opened) at the specified time t=t2, and the gas flow from the damaged pipeline should be redirected to the parallel pipeline's undamaged pipeline. Determining the time t=t2 becomes a critical issue at this point.

In the second part ($\ell_1 \leq x \leq \ell_3$), a venting process occurs at the point of rupture of the pipeline's integrity (x=$\ell_2$) due to the gas escaping into the atmosphere (Fig. 4).

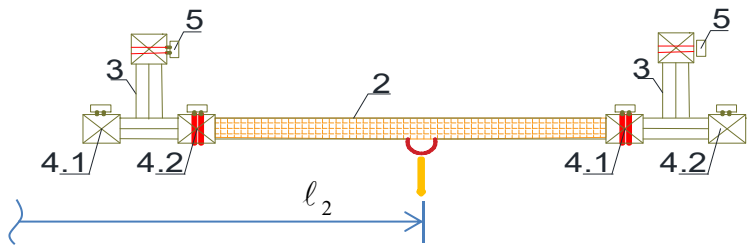

**Figure 4.** The principal scheme of the technological operation mode in the second part of the complex gas pipeline in the accident mode.

In this case, it is not possible to prevent the loss of the escaped gas. This is because both valves (sleeves) at the starting and ending points of that part (4.2 and 4.1) are closed. Consequently, due to the closure of these sleeves, the damaged section separates from the main part of the pipeline. At the same time, repair (welding) of the section with a damaged seal involves hazardous gas operations, so it must be carried out in accordance with the relevant regulations and instructions.

In the 3rd section of the pipeline ($\ell_3 \leq x \leq L$), the condensation process occurs due to gas consumption by consumers (Fig. 5). As seen in the scheme, the closure of the valve at the x= $\ell_3$ point of the damaged pipeline and the valves on the connections -3 result in a decrease in the pressure along the section over time, leading to a reduction in the gas stream's purity.

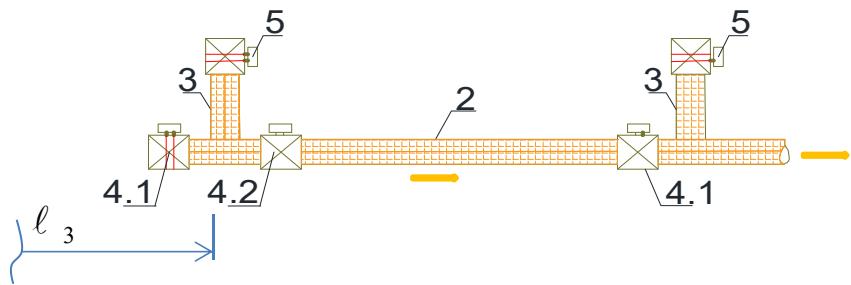

**Figure 5.** Principal Scheme of the Technological Operating regime in the 3rd section of the Complex Gas Pipeline in emergency mode.

The decrease in gas consumption at the end of the pipeline leads to a decrease in the reliability of gas supply to consumers. In other words, consumers cannot be supplied with the necessary gas pressure. After a certain time, the gas supply regime for consumers can be completely disrupted. At this point, the drawers (taps) at x= $\ell_1$ and x= $\ell_3$ should be activated at the -3 position, and the automatic valves on the -5 should be opened, allowing the required portion of gas from the undamaged pipe belt to flow into the undamaged third part of the damaged pipe where the accident occurred. Thus, consumers in the pipeline are continuously supplied with gas. At t = $t_2$, it is essential to activate the drawers (taps) installed on the connectors.

To assess the feasibility of the proposed scheme in the event of an accident, it is necessary to determine its dynamic state in all processes according to the operating principle. In other words, providing a mathematical solution for the upcoming physical processes involves determining the pressure of gas along the pipeline over time. For the mathematical solution of the processes, we use the following transfer functions [1, 2, 6].

For 1st part ($0 \leq x \leq \ell_1$).

$$\frac{\partial^2 \tilde{P}_1}{\partial X^2} = \frac{2a}{c^2}\frac{\partial \tilde{P}_1}{\partial t}$$

(1)

For 2nd part ($\ell_1 \leq x \leq \ell_3$).

$$\frac{\partial^2 \tilde{P}_2}{\partial X^2} = \frac{2a}{c^2}\frac{\partial \tilde{P}_2}{\partial t} + 2a\tilde{G}_{ut}(t)\delta(x-\ell_2)$$

(2)

For 3rd part ($\ell_3 \leq x \leq L$).

$$\frac{\partial^2 \tilde{P}_3}{\partial X^2} = \frac{2a}{c^2}\frac{\partial \tilde{P}_3}{\partial t}$$

(3)

Here, $\tilde{G} = \rho v$; $\tilde{P} = \rho \cdot c^2$

The mass flow measurements at the leakage point at various times represent the values of mass expenditures for the failure mode of the $\tilde{G}_s(t)$ pipeline belt, $\frac{Pa \ sec}{m}$.

$\delta(x-\ell_2)$ - Dirac's Function

c- For an isothermal process, the speed of sound propagation of gas is expressed in m/s

$\tilde{P}$ - Pressure at the narrowest points of the gas pipeline, Pa.

x-coordinate of the gas pipeline along its length, in kilometers.

t - time coordinate, seconds.

- hydraulic resistance coefficient.

- average density of gas, kg/m³.

v - average velocity of the gas flow, m/s.

d - diameter of the gas pipeline, meters.

The distribution of the function at the initial conditions for the physical process under consideration is given at the starting time t = 0. Since the process we are examining is a transient process, we will consider the connection of automatic valves installed at the -4.2 and -4.1 points on both sides of the x= $\ell_2$ point as the initial condition. If these valves are connected at time $t=t_1$, then the distribution of the sought function (pressure) at $t=t_1$ is given without taking into account inertial forces.

For $t=t_1$; $\tilde{P}_i(x,t_1) = P_i(x,t_1) \rightarrow i = 1, 2, 3$

The boundary conditions provide the initially known law-compliant variations in expenditure at four points along the pipeline during the considered time of the process. Proper specification of the boundary conditions completes the mathematical model of the process and enables a comprehensive and precise investigation of the occurring physical events.

at the point of $\dfrac{\partial \tilde{P}_1(x,t_1)}{\partial x} = -2a\tilde{G}_0(t)$ ; at the point of $\dfrac{\partial \tilde{P}_3(x,t_1)}{\partial x} = -2a\tilde{G}_s(t)$

at the point of $\dfrac{\partial \tilde{P}_1(x,t_1)}{\partial x} = 0$ ; at the point of $\dfrac{\partial \tilde{P}_2(x,t_1)}{\partial x} = 0$

$\dfrac{\partial \tilde{P}_2(x,t_1)}{\partial x} = 0$ ; $\dfrac{\partial \tilde{P}_3(x,t_1)}{\partial x} = 0$

Here, $\tilde{G}_0(t)$ and $\tilde{G}_s(t)$ represent the values of mass flow measured at different times at the initial and final points of the gas pipeline, $\dfrac{Pa \ \sec}{m}$.

When the stability of gas pipelines is compromised in experiments, preference is given to the Laplace transformation method. When referring to Laplace transformation $P(x,S) = \int_0^\infty P(x,t)e^{-st}dt$. $S=a+ib$ represents a complex number, where, $i=\sqrt{-1}$ is an imaginary number.

It is clear that through the application of Laplace transformation, equations (1), (2), and (3) are transformed into second-order differential equations, and their general solutions will be as follows.

$$\tilde{P}_1(x,s) = \dfrac{P_1(x,t_1)}{S} + c_1 Sh\lambda x + c_2 Ch\lambda x \qquad 0 \leq x \leq \ell_1$$

(4)

$$\tilde{P}_2(x,S) = \frac{P_2(x,t_1)}{S} + c_3 Sh\lambda x + c_4 Ch\lambda x - \beta \tilde{G}_{ut}(S) \int_{\ell_1}^{x} \delta(y-\ell_2) sh\lambda(y-x) dy$$

$\ell_1 \leq x \leq \ell_3$ (5)

$$\tilde{P}_3(x,S) = \frac{P_3(x,t_1)}{S} + c_5 Sh\lambda x + c_6 Ch\lambda x \quad \ell_3 \leq x \leq L$$

(6)

Here, $\lambda = \sqrt{\dfrac{2as}{c^2}}$, $\beta = \sqrt{\dfrac{2ac^2}{s}}$

Applying the Laplace transformation to the initial and boundary conditions, considering equations (4), (5), and (6), we can determine the constants C1, C2, C3, C4, C5, and C6. By substituting the values of these constants into equations (4), (5), and (6), we obtain the transformed form of the dynamic behavior of the considered process, in other words, the Laplace transform of the pressure distribution along the pipeline for the studied sections.

$$\tilde{P}_1(x,S) = \frac{P_1(x,t_1)}{S} + \beta \tilde{G}_0(S) \frac{ch\lambda(x-\ell_1)}{sh\lambda\ell_1} + \frac{dP_1(0,t_1)}{dx} \frac{ch\lambda(x-\ell_1)}{S\lambda sh\lambda\ell_1}$$

$0 \leq x \leq \ell_1$ (7)

$$\tilde{P}_2(x,S) = \frac{P_2(x,t_1)}{S} - \beta \tilde{G}_{ut}(S) \frac{ch\lambda(\ell_2-\ell_3)ch\lambda(x-\ell_1)}{sh\lambda(\ell_3-\ell_1)} -$$

$$- \begin{cases} 0 \to \ell_1 \leq x \leq \ell_2 \\ \beta \tilde{G}_{ut}(S) sh\lambda(x-\ell_2) \to \ell_2 \leq x \leq \ell_3 \end{cases}$$

(8)

$$\tilde{P}_3(x,S) = \frac{P_3(x,t_1)}{S} - \beta \tilde{G}_s(S) \frac{ch\lambda(x-\ell_3)}{sh\lambda(L-\ell_3)} - \frac{dP_3(L,t_1)}{dx} \frac{ch\lambda(x-\ell_3)}{S\lambda sh\lambda(L-\ell_3)}$$

$\ell_3 \leq x \leq L$ (9)

Based on the solution of the considered problem, since it is more convenient to determine the mathematical expression for the non-stationary flow of the gas stream in each of the three sections, we use the reverse Laplace transformation procedure for equations (7), (8), and (9) to find the original solution for each boundary. Thus, according to the initial and boundary conditions, we obtain the mathematical expressions for the distribution of pressure along the pipeline over time.

$$\tilde{P}_1(x,t) = P_1(x,t_1) + \frac{c^2}{\ell_1}\int_{t_1}^{t}\tilde{G}_0(\tau)d\tau + \frac{2c^2}{\ell_1}\sum_{n=1}^{\infty}\cos\frac{\pi n x}{\ell_1}\int_{t_1}^{t}\tilde{G}_0(\tau)\cdot e^{-\alpha_3(t-\tau)}d\tau +$$

$$+ \frac{c^2}{2a\ell_1}\frac{dP_1(0,t_1)}{dx}\int_{t_1}^{t}\left[1 + 2\sum_{n=1}^{\infty}\cos\frac{\pi n x}{\ell_1}e^{-\alpha_3\tau}\right]d\tau$$

$0 \leq x \leq \ell_1$ (10)

$$\tilde{P}_2(x,t) = P_2(x,t_1) - \frac{c^2}{(\ell_3-\ell_1)}\int_{t_1}^{t}\tilde{G}_{ut}(\tau)\left[1 + 2\sum_{n=1}^{\infty}(-1)^n\cos\frac{\pi n(\ell_2-\ell_3)}{\ell_3-\ell_1}\cdot\cos\frac{\pi n(x-\ell_1)}{\ell_3-\ell_1}e^{-\alpha_4(t-\tau)}d\tau\right]$$

$\ell_1 \leq x \leq \ell_3$ (11)

$$\tilde{P}_3(x,t) = P_3(x,t_1) - \frac{c^2}{(L-\ell_3)}\int_{t_1}^{t}\tilde{G}_s(\tau)\left[1 + 2\sum_{n=1}^{\infty}(-1)^n\cos\frac{\pi n(x-\ell_3)}{L-\ell_3}\cdot e^{-\alpha_5(t-\tau)}\right]d\tau -$$

$$- \frac{c^2}{2a(L-\ell_3)}\frac{dP_3(L,t_1)}{dx}\int_{t_1}^{t}\left[1 + 2\sum_{n=1}^{\infty}(-1)^n\cos\frac{\pi n(x-\ell_3)}{L-\ell_3}\cdot e^{-\alpha_5\tau}\right]d\tau$$

$\ell_3 \leq x \leq L$ (12)

Here, $\alpha_3 = \frac{\pi^2 n^2 c^2}{2a\ell_1^2}$, $\alpha_4 = \frac{\pi^2 n^2 c^2}{2a(\ell_3-\ell_1)^2}$, $\alpha_5 = \frac{\pi^2 n^2 c^2}{2a(L-\ell_3)^2}$

For t=t$_2$ duration, it is possible to assume $\tilde{G}_0(\tau) = \tilde{G}_0(\tau) =$ constant. In this case, we accept the laws expressed in equations (10), (11), and (12) as well as the given conditions to determine the regularity of the change in gas pressure for each of the three different sections of the complex pipeline's technological process.

P$_b$ = 14×10$^{-2}$ MPa ; P$_s$ =11 × 10$^{-2}$ MPa ; G$_0$= 10 $\frac{Pa\ sec}{m}$ ;

$2a = 0{,}1\frac{1}{\sec}$ ; $c = 383{,}3\frac{m}{\sec}$ ; L =3 ×10$^4$ m; $\ell_1$ = 1× 10$^4$ m ; $\ell_2$ = 1,45× 10$^4$ m; $\ell_3$ = 2× 10$^4$ m ;

Firstly, assuming the location of the gas pipeline leakage point is $\ell_2$ = 1,45 ×10$^4$ m and the closing time of automatic valves is t=t$_1$= 300 seconds, we accept (22) and (23) equations. By using these equations, we calculate values for $P_1(x,t_1), P_2(x,t_1)$ and $P_3(x,t_1)$ at every 5000 m and note them in the table below.

| X, km | 0 | 5 | 10 | 14,5 | 20 | 25 | 30 |
|---|---|---|---|---|---|---|---|
| $P_1(x,300) \cdot 10^{-2}$ MPa | 13,36 | 12,82 | 12,19 | 11,56 | 11,24 | 10,86 | 10,40 |
| $P_2(x,300), 10^{-2}$ MPa | - | - | 12,19 | 11,56 | 11,24 | - | - |
| $P_3(x,300), 10^{-2}$ MPa | - | - | - | - | 11,24 | 10,86 | 10,40 |

Using the data from Table 1, we calculate the values of the $\tilde{P}_1(x,t)$ expression for every 60 seconds and every 5 km for the first section of the complex gas pipeline, noting them in the table below, Table 2.

| x, km | $\tilde{P}_1(x,t)$, $10^{-2}$ Mpa | | | | | | | | | | |
|---|---|---|---|---|---|---|---|---|---|---|---|
| | t=0 | t=60 | t=120 | t=180 | t=240 | t=300 | t=360 | t=420 | t=480 | t=540 | t=600 |
| 0 | 13,36 | 14,13 | 14,58 | 15,02 | 15,46 | 15,91 | 16,35 | 16,79 | 17,24 | 17,68 | 18,13 |
| 5 | 12,82 | 13,22 | 13,67 | 14,11 | 14,55 | 15,0 | 15,44 | 15,89 | 16,33 | 16,77 | 17,22 |
| 10 | 12,19 | 12,47 | 12,91 | 13,36 | 13,8 | 14,24 | 14,69 | 15,13 | 15,57 | 16,02 | 16,46 |

Using Table 2, we construct the graph of the pressure distribution over time for the starting (x=0), middle (x=5 km), and end (x=10 km) points of the first section of the complex gas pipeline (Graph 1).

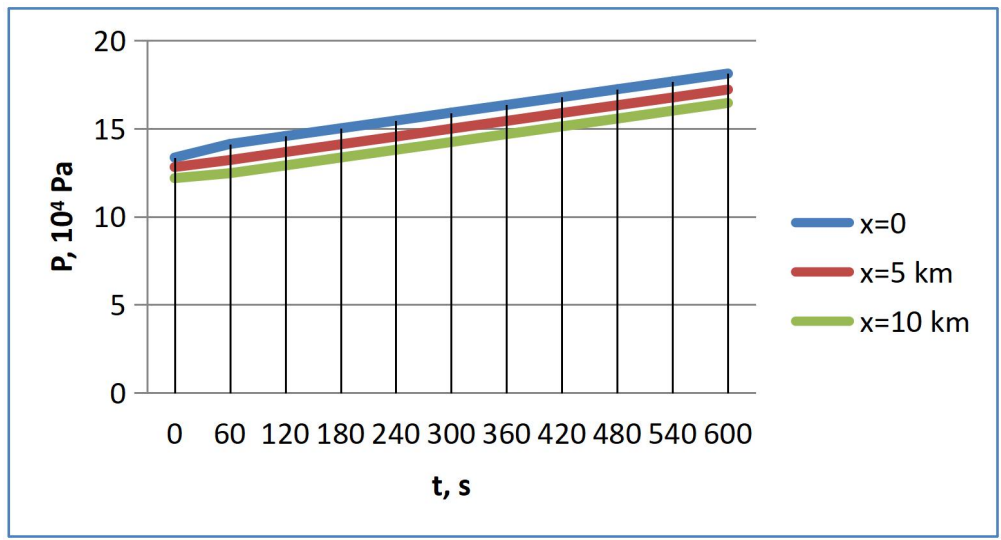

From Graph 1, it is evident that at time $t=t_1$, after the closing of the valves located at x= $\ell_1$ and x= $\ell_3$ points of the compromised pipeline, the pressure along the line

x= $\ell_1$ experiences an increase over time. The rate of increase in pressure is calculated at these coordinates. For example, at the starting point x=0, the recorded pressure value is 14.13 × 10⁴ Pa at t=60 seconds. This value increases to 15.91 × 10⁴ Pa at t=300 seconds and further to 18.13 × 10⁴Pa at t=600 seconds. Thus, the ratio of the increase at t=600 seconds is 1.28. At x=5 km, the pressure recorded at t=60 seconds is 13.22 × 10⁴Pa, and it increases to 17.22 × 10⁴Pa at t=600 seconds, resulting in a ratio of 1.3. Similarly, at x=10 km, the pressure increases from 12.47 × 10⁴ Pa at t=60 seconds to 16.46 × 10⁴Pa at t=600 seconds, with a ratio of increase being 1.32. Based on the analysis, it can be observed that the ratio of pressure increase along the axis from the starting point to the endpoint of the first section is increasing.

As mentioned, the compressors used at the beginning of the gas pipeline must be activated at points x= $\ell_1$ and x= $\ell_3$ to protect against the impact of pressure increase. This activation should occur precisely at t=t2, meaning the duration of this time should be determined and justified theoretically.

Continuously, using the data from Table 1, we calculate the values of the $\tilde{P}_2(x,t)$ expression for the starting point (x=10 km), the leakage point (x=14.5 km), and the endpoint (x=20 km) of the second section of the complex gas pipeline every 60 seconds, and record them in Table 3.

| x, km | $\tilde{P}_2(x,t)$, 10⁻² Mpa | | | | | | | | | | |
|---|---|---|---|---|---|---|---|---|---|---|---|
| | t=0 | t=60 | t=120 | t=180 | t=240 | t=300 | t=360 | t=420 | t=480 | t=540 | t=600 |
| 10 | 12,19 | 11,77 | 11,32 | 10,87 | 10,43 | 9,98 | 9,54 | 9,1 | 8,65 | 8,21 | 7,77 |
| 14,5 | 11,56 | 11,03 | 10,59 | 10,15 | 9,7 | 9,26 | 8,81 | 8,37 | 7,93 | 7,48 | 7,04 |
| 20 | 11,24 | 10,86 | 10,42 | 9,97 | 9,53 | 9,09 | 8,64 | 8,2 | 7,75 | 7,31 | 6,87 |

Using Table 3, we create the graph of the distribution of pressure over time for the starting point (x=10 km), leakage point (x=14.5 km), and endpoint (x=20 km) of the second section of the complex gas pipeline (Graph 2).

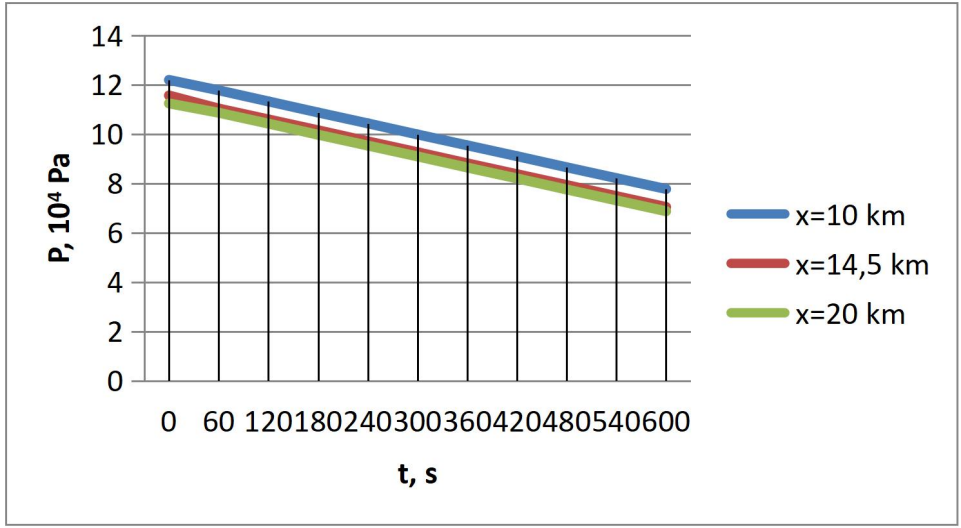

From Graph 2, it is evident that at t=$t_1$, after the closing of the control valves located at x= $\ell_1$ and x= $\ell_3$, the pressure along the line decreases over time for the second section of the gas pipeline ($\ell_1 \leq x \leq \ell_3$)The rate of this decrease is higher at the point of valve closure. In this case, at the starting point x=10 km, the pressure, noted at t=60 seconds, decreases from 11.77·10^4 Pa to 9.98·10$^4$Pa at t=300 seconds, and further decreases to 7.77·10$^4$Pa at t=600 seconds. In other words, the reduction ratio at t=600 seconds is 1.51 for x=14.5 km. At t=60 seconds, the recorded pressure value is 11.03×10$^4$ Pa, and at t=600 seconds, the pressure value decreases to 7.04×10$^4$ Pa. In other words, the reduction ratio is 1.57 for x=20 km. At t=60 seconds, the recorded pressure value is 10.86×10$^4$ Pa, and at t=600 seconds, the pressure value decreases to 6.87×10$^4$ Pa. It means the reduction ratio is 1.58.. Analysis suggests that the reduction ratio of pressure increases in the direction from the starting point to the final point of the second section.

As noted, it is impossible to release gas into the environment as a result of leakage. To determine the mass loss of the gas due to leakage, we perform the following consecutive operations based on equation (7):

$$\tilde{P}_2(x,S) - P_2(x,t_1) - \frac{c^2}{\ell_3 - \ell_1} \tilde{G}_{ut}(S) \frac{ch\sqrt{\frac{2as}{c^2}}(\ell_2 - \ell_3) ch\sqrt{\frac{2as}{c^2}}(x - \ell_1)}{\frac{sh\sqrt{\frac{2as}{c^2}}(\ell_3 - \ell_1)}{\sqrt{\frac{2as}{c^2}}(\ell_3 - \ell_1)}}$$

(13)

Here, $\lambda = \sqrt{\dfrac{2as}{c^2}}$

If we look at the expression in the (13) statement in terms of $S \to 0$ (expressed in the transformed form), then we adopt the following rule.

$$\int_{t_1}^{\infty} \frac{\partial \tilde{P}_2(x,t)}{\partial t} dt = -\frac{c^2}{\ell_3 - \ell_1} \int_{t_1}^{\infty} \tilde{G}_{ut}(t) dt$$

Here,

$$\int_{t_1}^{\infty} \tilde{G}_{ut}(t) dt = -\frac{\ell_3 - \ell_1}{c^2} \int_{t_1}^{\infty} \frac{d\tilde{P}_2(\ell_3, t)}{dt} dt$$

or,

$$\int_{t_1}^{\infty} \tilde{G}_{ut}(t) dt = -\frac{\ell_3 - \ell_1}{c^2} \tilde{P}_2(\ell_3, t)$$

(14)

Therefore, this gas loss is inevitable. From equation (14), it can be seen that the amount of this loss depends on the distance between displacements. As the distance between them increases, the amount of lost gas loss will also increase.

In conclusion, using the data from Table 1, we calculate the values of the expression $\tilde{P}_3(x,t)$ at the starting point (x=20 km), the midpoint (x=25 km), and the endpoint (x=30 km) of the third segment of the complex gas pipeline, at intervals of 60 seconds. We record these values in the following.

| x, km | $\tilde{P}_3(x,t)$, $10^{-2}$ Mpa | | | | | | | | | | |
|---|---|---|---|---|---|---|---|---|---|---|---|
| | t=0 | t=60 | t=120 | t=180 | t=240 | t=300 | t=360 | t=420 | t=480 | t=540 | t=600 |
| 20 | 11,24 | 10,96 | 10,52 | 10,08 | 9,63 | 9,19 | 8,74 | 8,3 | 7,86 | 7,41 | 6,97 |
| 25 | 10,86 | 10,46 | 10,01 | 9,57 | 9,13 | 8,68 | 8,24 | 7,8 | 7,35 | 6,91 | 6,46 |
| 30 | 10,4 | 9,63 | 9,19 | 8,74 | 8,3 | 7,85 | 7,41 | 6,97 | 6,52 | 6,08 | 5,63 |

Using Table 4, we create a graph depicting the distribution of pressure over time at the starting point (x=20 km), midpoint (x=25 km), and endpoint (x=30 km) of the third segment of the complex gas pipeline.

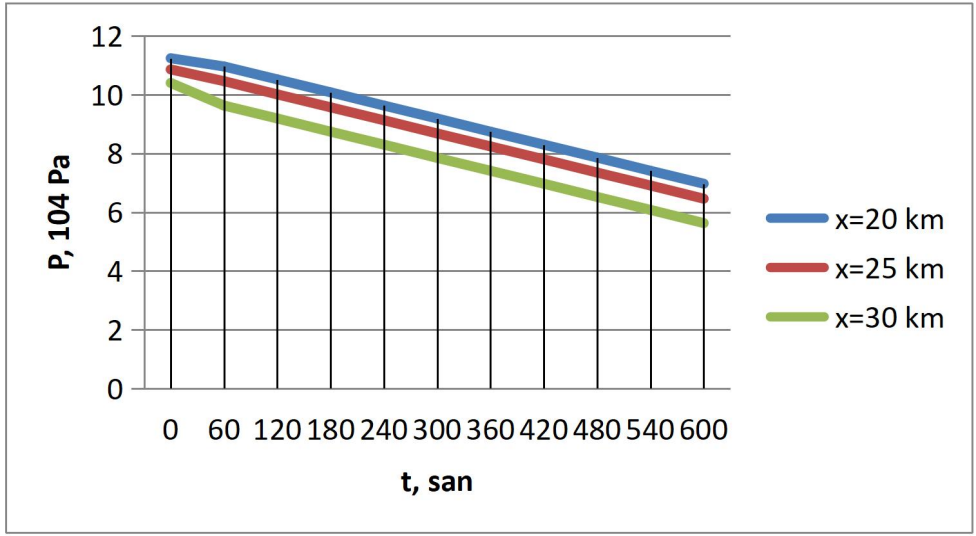

From Graph 3, it is apparent that at time $t=t_1$, after connecting the automatic valves located at points $x=\ell_1$ and $x=\ell_3$ on the depressurized pipeline, the pressure along the length of the gas pipeline's third segment ($\ell_3 \leq x \leq L$) decreases over time.

The rate of this decrease is calculated at specific coordinates. Specifically, at the starting point x=20 km, the pressure value decreases from $10.96 \cdot 10^4$ Pa at t=60 seconds to $9.19 \cdot 10^4$ Pa at t=300 seconds, and further decreases to $6.97 \cdot 10^4$ Pa at t=600 seconds. Exactly, at the point x=25 km, the pressure decreases from $9.63 \cdot 10^4$ Pa at t=60 seconds to $6.46 \cdot 10^4$ Pa at t=600 seconds, indicating a relative decrease of 1.57 at t=600 seconds. Indeed, at the endpoint x=30 km, the pressure decreases from $10.86 \cdot 10^4$ Pa at t=60 seconds to $5.63 \cdot 10^4$ Pa at t=600 seconds, resulting in a relative decrease of 1.93. Based on the analysis, we can observe that the relative decrease in pressure increases from the starting point to the endpoint in the direction of the third segment. As noted, during this time, the volume of gas supplied to consumers also decreases. In other words, after t=600 seconds, the volume of gas supplied to consumers will decrease by an average of 48 percent.

In order to improve the gas supply to consumers and restore its continuity, it is essential to activate the automatic valves installed on connectors at points $x=\ell_1$ and $x=\ell_3$ at the specified time $t=t_2$. The determination and justification of this duration rely on utilizing expression (10) to identify the lawfulness of the pressure change at the starting point of the gas pipeline.

$$\tilde{P}_1(0,t) = P_1(0,t_1) + \frac{2c^2\tilde{G}_0}{\ell_1}(t-t_1) + \frac{2a\ell_1}{3}\tilde{G}_0 - \frac{4a\ell_1}{\pi^2}\tilde{G}_0 \sum_{n=1}^{\infty} \frac{e^{-\alpha_3(t-t_1)} + e^{-\alpha_3 t} - e^{-\alpha_3 t_1}}{n^2}$$

(15)

The simplified expression of formula (15) is as follows.

$$\tilde{P}_1(0,t) = P_1(0,t_1) + \frac{2c^2\tilde{G}_0}{\ell_1}(t-t_1) + \frac{2a\ell_1}{3}\tilde{G}_0 -$$

$$- \frac{4a\ell_1}{\pi^2}\tilde{G}_0 \left[ e^{-\alpha_3(t-t_1)} + e^{-\alpha_3 t} - e^{-\alpha_3 t_1} - \alpha(t-t_1)(1-C) \right]$$

(16)

Here, "C" is a Euler's formula, and its value is taken as $C = 0,577215$.

Analyses of the processed technical and economic parameters of gas compressor units have shown that when ε approaches 1.35 (ε ⟩ 1,35), there is only one optimal value satisfying this condition

As a result, it is clear that the basis of the normal operating mode of the compressor station is the ratio ε, which reflects the number of gas compressor units and their compression degree. Generally, since the main equipment in the exploitation of complex main gas pipelines is the compressor station, significant attention should be paid to the limited value of ε for the protection of their crucial compressor units. As a conclusion of the analyses conducted above, it can be stated that the theoretical determination of the duration t=t₂ should be based on the ε ratio. In this case, for the determination of t₂, the condition represented by the change in the initial pressure of the gas pipeline we analyzed in stationary mode - $P_b$ and the gas pipeline's first section in the accident mode (after the automatic valves are connected) denoted by $\tilde{P}_1(o,t)$ should be taken into account.

$$\varepsilon \geq \frac{\tilde{P}_1(o,t)}{P_b}$$

(17)

If we replace the "≥" symbol with "=", we obtain the following equality for t= t₂ in expression (17).

$$\tilde{P}_1(o,t) = \varepsilon P_b$$

In that case, considering expression (16) in relation to expression (17), we derive the following analytical rule for determining the time of activating the automatic valves installed on connectors at points x=$\ell_1$ and x=$\ell_3$.

$$t_2 = t_1 + \frac{\ell_1}{\tilde{G}_0} \frac{\varepsilon P_b - P(0,t_1) - 2a\ell_1 \left(\frac{1}{3} - \frac{1}{\pi^2}\right)}{(2-C)}$$

(18)

Since the process we are examining is a transition process, we have assumed the initial condition as t= $t_1$. It is evident that, when the valve is opened at point x= $\ell_2$ in the complex gas pipeline, the equations expressing the distribution of pressure for all three segments until the main connection of automatic valves at the left -4.2 and right -4.1 sections are as follows [5].

$$P_{1,2}(x,t) = P_1 - 2aG_0 x + \frac{8aL}{\pi^2} G_o \sum_{n=1}^{\infty} \frac{e^{-\alpha_1 (2n-1)^2 t}}{[(2n-1)]^2} \cos\frac{\pi(2n-1)x}{L} -$$

$$\frac{4aL}{\pi^2} G_o \sum_{n=1}^{\infty} \left[1-(-1)^n\right] \frac{e^{-\alpha_1 n^2 t}}{n^2} \cos\frac{\pi n x}{L} -$$

$$-\frac{c^2 t_1}{L} G_{ut} - 2aG_{ut} \left[\frac{x^2}{2L} + \frac{\ell_2^2}{2L} + \frac{L}{3} - \ell_2\right] + \frac{4aL}{\pi^2} G_{ut} \sum_{n=1}^{\infty} \cos\frac{\pi n x}{L} \cos\frac{\pi n \ell_2}{L} \cdot e^{-\alpha_1 n^2 t}$$

(19)

The solution of the expressions for the infinite path in rule (21) is as follows.

$$\sum_{n=1}^{\infty} \left[\frac{e^{-\alpha_1 (2n-1)^2 t}}{[(2n-1)]^2} \cos\frac{\pi(2n-1)x}{L} - \left[1-(-1)^n\right]\frac{e^{-\alpha_1 n^2 t}}{n^2} \cos\frac{\pi n x}{L}\right] = 0$$

If we consider this simplification in expression (21), it becomes as follows.

$$P_{1,2}(x,t) = P_1 - 2aG_0 x - \frac{c^2 t_1}{L} G_{ut} - 2aG_{ut} \left[\frac{x^2}{2L} + \frac{\ell_2^2}{2L} + \frac{L}{3} - \ell_2\right] +$$

$$+\frac{4aL}{\pi^2} G_{ut} \sum_{n=1}^{\infty} \cos\frac{\pi n x}{L} \cos\frac{\pi n \ell_2}{L} \cdot \frac{e^{-\alpha_1 n^2 t}}{n^2}$$

(20)

$$P_3(x,t) = P_1(x,t) - 2aG_{ut}\ell_2 + 2aG_{ut}x$$

(21)

Here, $\alpha_1 = \dfrac{\pi^2 nc^2}{2aL^2}$

By using expression (20), we can determine the value of the expression. Based on Table 1, the value of the expression at t=300 seconds is =13.36·0⁴Pa.
In that case, considering ε=1.35 and using the provided values and expression (18), we can determine the time for activating the automatic valves installed on connectors.

$$t_2 = t_1 + \frac{\ell_1}{c^2 \tilde{G}_0} \frac{\varepsilon P_b - P(0,t_1) - 2a\ell_1 \tilde{G}_0 \left(\frac{1}{3} - \frac{1}{\pi^2}\right)}{(2-C)} = 300 + 255 = 555 \text{ sec}$$

So, if the proposed scheme and the given parameters are used, and if the valves installed on the damaged pipeline are closed at $t_1 = 300$ seconds when the gas pipeline is depressurized, then the valves on the connectors should be opened at $t_2 = 255$ seconds. If we consider the moment of the accident, this duration t2 will be 555 seconds.

**Conclusion**

Taking into account the transition duration of the motion equation based on classical methods, the solution of the system of differential equations describing the non-stationary gas flow in a pipeline with a constant minimum diameter has been achieved. The solution, introduced at the initial time of the solution, allows meeting the boundary conditions' requirement of zero equality, enabling obtaining the comprehensive notation of the analytical model. However, it does not limit the scope of the model for sharp changes in gas flow or pressure.

The methodology for calculating transient processes described both by arbitrary boundary conditions (closing of valves) and heat transfer equations in main gas pipelines has been provided. The method for calculating changes in pressure along the pipeline for various operating modes of the studied section of a complex gas pipeline has been developed and prepared.

The proposed methods can be widely used in the design, reconstruction, and creation of automated technological control systems for main gas pipelines. An analytical approach has been employed to determine the time of transition from one process to another in complex gas pipelines. The adopted principle is suitable for the application of automated control systems for technological processes in pipelines and is convenient for engineering calculations.

Examples of calculations and their comparison with theoretical sources allow evaluating the high accuracy of the proposed methods for calculating transient processes in complex main gas pipelines.

**REFERENCES**

[1]. A. G. Vanchin, Methods for calculating the operation mode of complex main gas pipelines / A. G. Vanchin // Oil and Gas Business. 2014. No. 4. P. 192–214.

[2].S. Trofimov, V. A. Vasilenko, E. V. Kocharyan, Quasilinearization of the gas flow equation in a pipeline, Oil and Gas Business, no. 1, 1–11, 2003

[3].A.A. Apostolov, A.S. Verbilo, B.S. Pankratov, Automation of dispatch control of a gas transportation enterprise, Overview. Ser. Automation, Telemechanics, and Communication in the Gas Industry, M.: LLC "IRC Gazprom", 72, 1999

[4].S. Elaoud, B. Abdulhay, E. Hadj-Taieb, Archives of Mechanics, 66, 4, 269–286, 2014

[5].S. Fikov, The best estimate of the linearization parameter of the mathematical model of non-stationary gas flow in pipelines, Innovation. Education. Energy efficiency: materials of the XIV International Scientific and Practical Conference / GIPK "GAZ-INSTITUTE". Minsk, October 29–30, 82–85, 2020

[6].V. I. Panferov, S. V. Panferov, Modeling of non-stationary processes in gas pipelines / Bulletin of the South Ural State University. Series Construction and Architecture. 2007. No. 4. Issue 14. P. 44–47.

[7].I.Q. Aliyev, M.Z. Yusifov, N.I. Alizade, Investigation of a new method for determining the damage location in the analysis of non-stationary flow parameters of complex gas pipelines, no. 1, January 20, 2024 ALMATY, KAZAKHSTAN, ISSN 2709-1201.